\documentclass{aims}

\usepackage{txfonts}
\def\typeofarticle{Type of article}
\def\currentvolume{xxxx}
\def\currentissue{xxxx}
\def\currentyear{xxxx}

\def\ppages{xxx--xxx}
\def\DOI{xxxx}
\def\Received{xxxx}
\def\Accepted{xxxx}
\def\Published{xxxx}

\newtheorem{theorem}{Theorem}[section]
\newtheorem{lemma}[theorem]{Lemma}
\newtheorem*{mainthm}{Main Theorem}

\theoremstyle{definition}
\newtheorem{definition}[theorem]{Definition}

\theoremstyle{theorem}
\newtheorem{remark}[theorem]{Remark}

\numberwithin{equation}{section}

\usepackage{float}

\usepackage[utf8]{inputenc}
\usepackage[T1]{fontenc}
\usepackage{graphicx}
\usepackage{mathrsfs}
\usepackage{dsfont}
\usepackage{amssymb,graphicx,xcolor}
\usepackage{amssymb,bbm,enumerate}
\usepackage{esint}
\usepackage{colortbl}
\usepackage{xcolor}
\usepackage[utf8]{inputenc}

\usepackage{amsmath,amssymb,amstext,amsthm}

\usepackage{tikz}
\usepackage{amssymb}
\usepackage{bbm}

\usepackage{amsmath,amssymb,amstext,amsthm}

\newcommand{\fm}{\dot{H}^{-1}}

\DeclareMathOperator{\Lip}{Lip}
\DeclareMathOperator{\dist}{dist}

\DeclareMathOperator{\mix}{mix}

\def\Xint#1{\mathchoice
	{\XXint\displaystyle\textstyle{#1}}%
	{\XXint\textstyle\scriptstyle{#1}}%
	{\XXint\scriptstyle\scriptscriptstyle{#1}}%
	{\XXint\scriptscriptstyle\scriptscriptstyle{#1}}%
	\!\int}
\def\XXint#1#2#3{{\setbox0=\hbox{$#1{#2#3}{\int}$ }
		\vcenter{\hbox{$#2#3$ }}\kern-.6\wd0}}

\def\dashint{\Xint-}

\numberwithin{equation}{section}

\makeatletter
\renewcommand{\@biblabel}[1]{#1\hfill \hspace{-0.2cm}}
\makeatother

\begin{document}

\title{Sub-exponential mixing of generalized cellular flows \\ with bounded palenstrophy}

\author{%
  Gianluca Crippa\corrauth\affil{1} and
  Christian Schulze\affil{1}
}

% \shortauthors is used in copyright information in the end of the paper
\shortauthors{the Author(s)}

\address{%
  \addr{\affilnum{1}}{Department of Mathematics and Computer Science, University of Basel, 
  Spiegelgasse 1, 4051~Basel, Switzerland}
}

% corresponding author
\corraddr{gianluca.crippa@unibas.ch
}

\begin{abstract}
We study the mixing properties of a passive scalar advected by an incompressible flow. We consider a class of  cellular flows (more general than the class in [Crippa-Schulze M$^3$\!AS 2017]) and show that, under the constraint that the palenstrophy is bounded uniformly in time, the mixing scale of the passive scalar cannot decay exponentially.
\end{abstract}


\maketitle

\section{Introduction}

We are interested in the mixing properties of passive scalars advected by incompressible flows. How well a scalar is mixed by a flow is an important problem in fluid mechanics, with applications to atmospheric and oceanographic science, biology, and chemistry, just to name a few. In many situations the advected scalar has a negligible feedback on the underlying flow, and mathematically this results in the fact that the scalar solves a linear continuity equation with a given velocity field (see~\eqref{Cauchyprob}).

Let us now specify our mathematical setting. We consider a binary passive scalar $\rho\in\{+1,-1\}$ advected by a time-dependent, divergence-free velocity field $u$ on the two dimensional open unit square $\mathcal{Q}=(-\frac{1}{2},\frac{1}{2})^2\subset \mathbb{R}^2$, with $u=0$ on $\partial \mathcal{Q}$.
Given an initial condition $\bar{\rho}$, the scalar $\rho$ satisfies the Cauchy problem for the continuity equation with velocity field $u$, that is, it solves
\begin{equation}
\label{Cauchyprob}
\begin{cases} \partial_t\rho + \textnormal{div} (u\rho) =0 & \mbox{on } \mathbb{R}_{+}\times\mathbb{R}^2\\ \rho(0,\cdot)=\bar{\rho} & \mbox{on }\mathbb{R}^2.\end{cases}
\end{equation}
Outside of $\mathcal{Q}$, both the velocity field $u$ and the solution $\rho$ are identically zero for all times. 
Without loss of generality we can assume that $\bar{\rho}$ satisfies $\fint_{\mathcal{Q}}\bar{\rho}=0$ (otherwise we can simply consider $\tilde{\rho}=\rho-\fint\bar{\rho}$) and we observe that such a condition will be satisfied by the solution $\rho(t,\cdot)$ at all times.  The results presented in this paper hold in all dimensions $d\geq 2$ as well with the necessary changes in the scaling analysis. For the construction of examples, $d=2$ is the most restrictive case. 

\medskip

\noindent\textbf{Mixing scales.} In order to mathematically study the problem of optimal mixing we first need to specify how we measure the ``degree of mixedness'' of the solution~$\rho$ at each time.  The resulting ``mixing scales'' vary greatly in the literature and depend not only on the field of mathematics, but also on the specific PDE under consideration. The presence of diffusion, for instance, changes the mixing process drastically, and the $L^2$ norm is often a suitable mixing scale in the diffusive setting. However, the $L^2$ norm is conserved in the inviscid equation~\eqref{Cauchyprob} and therefore it is not a good measure of the mixing scale. One needs instead to find a quantification of the (possible) rate of weak convergence to zero of the solution. Two notions of mixing scale in such a case are available in the literature, namely the \textit{geometric} and the \textit{functional} mixing scale.

\begin{definition}[Geometric mixing scale~\cite{Bressan}]
	\label{Def:GeoMix}
	Given an accuracy parameter $0<\kappa<1$, the geometric mixing scale of $\rho(t,\cdot)$ is the infimum of all $\epsilon>0$ such that for every $x\in \mathbb{R}^2$ it holds
	\begin{equation}
	\frac{1}{\|\rho(t,\cdot)\|_{\infty}}\left|\fint_{B(x,\epsilon)}\rho(t,y)\,dy\,\right|<\kappa\,.
	\end{equation}
	Note that in the binary case~$\rho\in\left\{+1,-1\right\}$ the above condition is equivalent to
	\begin{equation}
	\label{AggroBerlin}
	\frac{1-\kappa}{2}< \frac{| \{\rho(t,\cdot) = +1\} \cap B(x,\epsilon)|}{|B(x,\epsilon)|}< \frac{1+\kappa}{2} \,.
	\end{equation}
	We denote by $\mathcal{G}(\rho(t,\cdot))$ such infimum. 
\end{definition}

\begin{definition}[Functional mixing scale~\cite{Mezic,Thieff,LinThif}]
	\label{Def:FuncMix}
	The functional mixing scale of $\rho(t,\cdot)$ is
	\begin{equation}
	\label{H1Duality}
	\|\rho(t,\cdot)\|_{\dot{H}^{-1}(\mathbb{R}^2)}=\sup\left\lbrace \int_{\mathbb{R}^2}\rho(t,x)\,\xi(x)\,dx \,:\, \|\nabla \xi \|_{L^2(\mathbb{R}^2)}\leq 1 \right\rbrace \, .
	\end{equation}
\end{definition}

\noindent\textbf{Notation.} It is often convenient to write~$\mix(\rho(t,\cdot))$ to denote any of the two mixing scales, in case a statement holds for both of them. 

\smallskip

\noindent It is worth mentioning that the two mixing scales are strongly related, but not equivalent, see the examples in~\cite{Lin}.

\medskip

\noindent\textbf{Optimal mixing.} The central question in optimal mixing is how fast the passive scalar can be mixed under suitable (energetic) constraints on the velocity field~$u$, that is, how fast (as a function of the time) the mixing scales can decay to zero. Answering this question entails two parts: on the one hand, one establishes lower bounds on the decay of the mixing scales, and on the other hand one exhibits (explicit) examples of velocity fields that mix optimally under such a constraint, that is, examples saturating the lower bounds.

To illustrate one particular, yet very relevant, example, let us focus on the case of a \textit{fixed enstrophy constraint}, that means $\|\nabla u(t,\cdot)\|_{L^2}\leq C$ for all $t\geq 0$. Under this constraint, Crippa and De Lellis~\cite{LellisCrippa} proved an exponential lower bound for the geometric mixing scale. This result was later extended by Iyer, Kiselev and Xu~\cite{Kiselev} and by Seis~\cite{Seis} to the functional mixing scale, so that we have
\begin{equation}
\label{EQ:LowerBounds}
\mathcal{G}(\rho(t,\cdot))\geq c\exp(-ct)
\hspace{0.5cm}\text{ and }\hspace{0.5cm}
\|\rho(t,\cdot)\|_{\fm}\geq c\exp(-ct)\,,
\end{equation}
where the constant $c$ depends on $\sup_t \|\nabla u(t,\cdot)\|_{L^2}$ and on the initial datum $\bar\rho$.

For future use, we record here that the proof of the lower bounds in~\eqref{EQ:LowerBounds} is based on a key lemma from~\cite{LellisCrippa} showing \textit{Lusin-Lipschitz regularity} of the flow of a Sobolev velocity field, in the form of a quantitative bound:
\begin{theorem}[Crippa and De Lellis]
	\label{Thm:CrippaDeLellisRegOG}
	Let $X:[0,T]\times\mathbb{R}^2\to \mathbb{R}^2$ be the flow associated to a velocity field in $L^1([0,T];W^{1,p}(\mathbb{R}^2))$, with $p>1$. Then, for every $\epsilon>0$ and every $R>0$, we can find a set $K\subset  B(0,R)$ such that $|B(0,R)\setminus K |\leq \epsilon$ and for any $0\leq t\leq T$ we have
	\begin{equation}
	\label{Lippi}
	\Lip(X(t,\cdot)|_{K})\leq \exp\left(\frac{c_p \, \int_0^t \|\nabla u(s,\cdot)\|_{L^p}\,ds }{\epsilon^{1/p}}\right)\,.
	\end{equation}
\end{theorem}

\noindent The above result somewhat extends the Cauchy-Lipschitz theory to the Sobolev framework. 

The  exponential lower bounds in~\eqref{EQ:LowerBounds} turn out to be sharp. Indeed, Alberti, Crippa and Mazzucato~\cite{AlbCrippa} and Yao and Zlato\v{s}~\cite{Yao} constructed explicit examples of entrophy-constrained velocity fields exhibiting exponential mixing. Both examples are of \textit{cellular type}, a special structure that plays an important role in our analysis and that we describe in the next paragraph. 

\medskip

\begin{figure}[h]
	\begin{center}
		\scalebox{0.5}{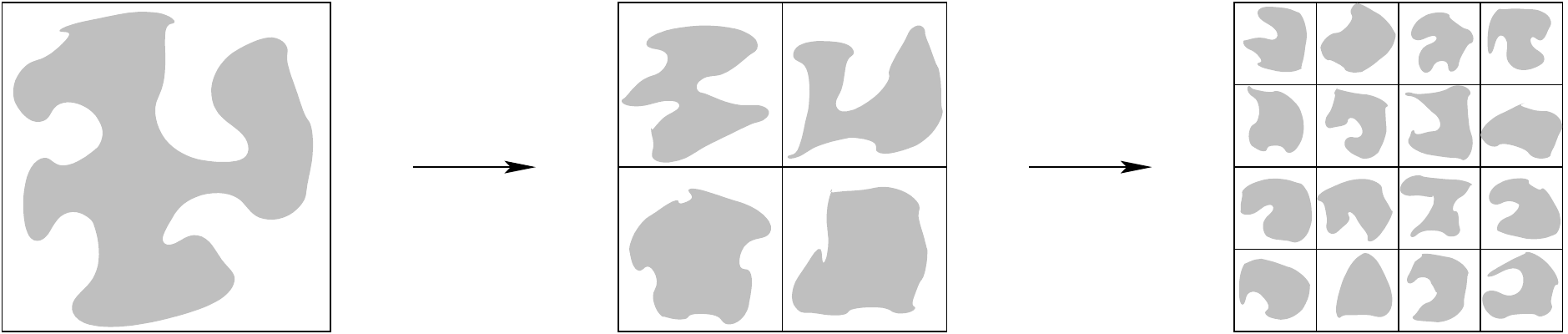}
	\end{center}
	\caption{First two steps of a cellular flow \label{concmus}}
\end{figure}

\noindent\textbf{Cellular structure.} The rough idea of the cellular structure is sketched in~Figure~\ref{concmus}. We start with a binary initial condition $\rho_0$ with zero average on the $2$-dimensional square $\mathcal{Q}$. In a first step we subdivide $\mathcal{Q}$ into $4$ disjoint sub-squares $D_1,\ldots,D_4$ of equal size. We then construct a velocity field $u_0$ in such a way that the solution $\rho(1,\cdot)$ has zero average on each of the sub-squares, which means that the tracer gets equally distributed among $D_1,\ldots,D_4$, as schematically visualized in~Figure~\ref{concmus}. 
We continue inductively by subdividing each $D_i$ into four equal sub-squares and distributing the tracer equally among those sub-squares with \textit{tracer movements localized in $D_i$} respectively, and so forth. In a sense, in each step the problem of mixing on a large area is transferred into the problem of \textit{separately} mixing on four smaller areas. 

Using a cellular type structure for the construction of mixing flows has many advantages. In addition to the convenient inductive nature of the construction, a crucial benefit of this structure is that it allows to keep track of both the geometric and the functional mixing scale at all times by simple scaling arguments, as we will see later in this paper. 

\medskip

\begin{figure}[H]
 	\begin{center}
 		\scalebox{0.6}{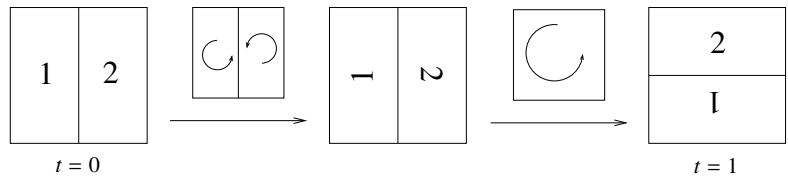}
 	\end{center}
 	\caption{Action of the velocity field $u$ by Elgindi and Zlato\v{s} for $0\leq t\leq 1$.\label{yyyy}}
 \end{figure}

\noindent\textbf{Fixed palenstrophy constraint.} In this paper we deal with a \textit{fixed palenstrophy constraint}, by which we mean a uniform-in-time bound on the $\dot{W}^{s,p}$ norm of the velocity field~$u$, where $s>1$ and $1< p\leq\infty$ (with a slight abuse of terminology, since more usually the term ``palenstrophy'' denotes the $\dot{H}^2$ norm of the velocity field $u$).  We refer to~\cite{Lin,LinErr} for more details on the concept of palenstrophy and its significance in the context of mixing of passive scalars. Under a fixed palenstrophy constraint we immediately inherit the exponential lower bounds~\eqref{EQ:LowerBounds} for both mixing scales from the fixed enstrophy case due to the (fractional) Poincar\'e inequality 
\begin{equation}
\label{Poincarefrac}
\|u(t,\cdot)\|_{\dot{W}^{s,p}(\mathbb{R}^2)}\geq C_{s,p} \|u(t,\cdot)\|_{\dot{W}^{1,p}(\mathbb{R}^2)}\, .
\end{equation}
 Elgindi and Zlato\v{s}~\cite{ElgindiZlatos} proved that this lower bound is sharp for a certain range of $s,p>1$. The authors prove that a variation of \textit{Baker's map} (a discrete map well-known in the theory of dynamical systems) is associated to a flow. The action of the velocity field~$u$ on the time interval~$[0,1]$ is sketched in Figure~\ref{yyyy}. The velocity field is then periodically extended and therefore it is \textit{not} of cellular type, since the same velocity field is repeated periodically in time without going to smaller scales.
\noindent
It is in fact believed that the exponential bound is sharp in the full range $s,p>1$, but no analytic examples are available (see also the numerical simulations in~\cite{LinThif}).

We further notice that it is possible to rescale (in time) the cellular flows in~\cite{AlbCrippa} and~\cite{Yao} in such a way that they satisfy the fixed palenstrophy constraint. However, by doing so, the rescaled flows mix at a slower, polynomial rate. It is therefore natural to ask whether such ``reduced mixing ability'' is specific to these particular examples (which are not only cellular, but in fact have a self-similar character), or if it applies to all cellular examples under this constraint and is due to the increased localization of the action of the flow.

We partially answered this question in~\cite{CellP}, where we considered cellular flows based on ``building blocks'' that do not mix ``too well'' at every iteration. To be more specific, consider again the sketch in Figure~\ref{concmus}. Note that, in every iteration, the side-length of the cells decays by a fixed factor~$1/2$. In~\cite{CellP} we considered cellular flows for which also the mixing scale decays \textit{at most} by a factor~$1/2$ in every iteration, as it is the case in Figure~\ref{concmus}. As a consequence, in the context of~\cite{CellP}, after every iteration the ratio between the mixing scale and the side length of the cells is constant. Note that the cellular examples in~\cite{AlbCrippa} and~\cite{Yao} follow this rule and therefore fall within the setting of~\cite{CellP}. Under this assumption, we proved that both mixing scales cannot decay faster than \textit{polynomially} under a fixed palenstrophy constraint.

The class considered in~\cite{CellP} is meaningful as it contains several explicit examples of mixers (for instance, self-similar and quasi self-similar ones, see again~\cite{AlbCrippa}). However, the requirement on the constant ratio between the cell size and the mixing scale is quite restrictive and it is not evident whether the sub-exponential lower bound in~\cite{CellP} is due to this condition or, more generally, it holds as a consequence of the increased localization of the velocity field. The natural questions that arise following this line of thinking are: What happens if we allow a cellular velocity field to mix arbitrarily well (with respect to the cell size) at every iteration? Is this restriction (weaker than the one in~\cite{CellP}) strong enough to exclude exponential mixing under fixed palentstrophy?

In the present article we simplify both the setup and the proof in~\cite{CellP} and completely drop the restriction preventing the building blocks from mixing ``too well'', still showing that exponential mixing cannot happen. We describe in detail the \textit{generalized cellular structure} we are able to deal with in Section~\ref{s:setting}, while in Section~\ref{s:main} we state and prove our Main Theorem.  In order to make the comparison with the statement in~\cite{CellP} more explicit, we underline that in~Definition~\ref{celltype}(ii) we assume that, after $n$ iterations, the solution is mixed at scale $\lambda^{n+\sigma(n)}$, where the increasing function $\sigma : {\mathbb N} \to {\mathbb N}$ may be unbounded. The result in~\cite{CellP} corresponds to the case $\sigma(n)=0$, as we comment in Remark~\ref{r:compa}(iii). The case of an unbounded function $\sigma$ corresponds to a mixing scale which, after $n$ steps, is strictly finer than $\lambda^n$.
 
\medskip

\noindent\textbf{Acknowledgements.} This research has been partially supported by the ERC Starting Grant 676675~FLIRT.

\section{Generalized Cellular Structure}
\label{s:setting}

Fix $\lambda>0$ such that $1/\lambda$ is an integer greater or equal than two. We define the \textit{tiling} of $\mathcal{Q} =(-\frac{1}{2},\frac{1}{2})^2$ corresponding to the parameter $\lambda$ as follows.
\begin{definition}
	\label{Def:tilingEX}
	We denote by $\mathcal{T}\!\!_\lambda$ the \textit{tiling} of $\mathcal{Q}$ with squares of side $\lambda$, consisting of the $1/\lambda^2$ open squares of the form 
	\begin{equation*}
	\left\lbrace (x,y)\in \mathcal{Q} \, :\, -\frac{1}{2}+k\lambda<x<-\frac{1}{2}+(k+1)\lambda \;\textnormal{ and }\; -\frac{1}{2}+h\lambda<y<-\frac{1}{2}+(h+1)\lambda\right \rbrace
	\end{equation*}
	for $k,h=0,\ldots,1/\lambda -1$.
\end{definition}

Let $(T_n)_{n\in\mathbb{N}}$ be an arbitrary (increasing) sequence of times. We define the generalized cellular structure in the following way. For times $T_n\leq t< T_{n+1}$, we assume that in any tile $Q\in \mathcal{T}\!\!_{\lambda^n}$ we can write the velocity field $u$ and corresponding solution $\rho$ as
\begin{equation*}
\label{TilVelSol}
u(t,x)=\frac{\lambda^n}{T_{n+1}-T_n}u_{Q,n}\left(\frac{t-T_n}{T_{n+1}-T_n},\frac{x-r_Q}{\lambda^n}\right)\hspace{0.5cm}\textnormal{and}\hspace{0.5cm}\rho(t,x)=\rho_{Q,n}\left(\frac{t-T_n}{T_{n+1}-T_n},\frac{x-r_Q}{\lambda^n}\right),
\end{equation*}
where $(u_{Q,n},\rho_{Q,n})$ solves~\eqref{Cauchyprob} for $0\leq t\leq 1$ and $r_Q$ denotes the center of $Q\in \mathcal{T}\!\!_{\lambda^n}$. Furthermore, we assume that $\rho_{Q,n}(0,\cdot)$ is both \textit{mixed} and \textit{un-mixed} at scale $\lambda^{\sigma(n)}$, where $(\sigma(n))_{n\in\mathbb{N}}\subset \mathbb{N}$ is an increasing sequence of positive integers. We  specify our notion of mixed and un-mixed in Definition~\ref{Def:MixedUnmixed}, after a short digression on some background material.

We adapt the notion of un-mixed from~\cite{CellP}. The idea is to call a binary function $\rho$ \textit{un-mixed at scale~$r$} if there exists a ball of radius~$r$ that violates the geometric mixing scale condition in~\eqref{AggroBerlin} significantly. For this we fix a parameter $\bar{\gamma}\in(0,1)$. This parameter will be fixed for the rest of this section and for this reason explicit dependencies of constants on $\bar{\gamma}$ will often not be stated.  In the following, $0<\kappa<1$ is an accuracy parameter as in Definition~\ref{Def:GeoMix}.
\begin{definition} [Characteristic length scale]
	\label{zulu}
	Let a set $J\subset E$ of positive measure be given. We denote by $\mathcal{C}=\mathcal{C}(E, J)$ the set of all admissible balls $B(x,r) \subset E$ such that
	\begin{equation}
	\label{saucinonu}
	\frac{|J\cap B(x,r)|}{|B(x,r)|}>1-\frac{1-\kappa}{2}\cdot\bar{\gamma}.
	\end{equation}
	We define the \textit{characteristic length scale} of the set $J$ with respect to $E$ as
	\begin{equation*}
	LS\!_{E}(J)=\sup\big\{ r>0 \,\text{ such that there exists }\, B(x,r)\in \mathcal{C}\big\}.
	\end{equation*}
\end{definition}
\begin{remark}
	This scale determines the largest radius of a ball which violates \emph{significantly} (that is, with the uniform gap $\bar{\gamma}$)
	condition \eqref{AggroBerlin} of the geometric mixing scale for a binary tracer $\rho\in \{+1,-1\}$ and will serve as our measurement of how un-mixed the tracer is. To see this, let us set~$\bar{\gamma}=1$ and
	\begin{equation}
	J=\{\rho(t,\cdot) = +1\}\,.
	\end{equation}
	 Note that \eqref{saucinonu} becomes
	 \begin{equation}
	\frac{| \{\rho(t,\cdot) = +1\} \cap B(x,r)|}{|B(x,r)|}> \frac{1+\kappa}{2}\,,
	 \end{equation}
	 and hence $B(x,r)$ violates condition~\eqref{AggroBerlin} of the geometric mixing scale. For~$\bar{\gamma}\in(0,1)$, it violates the condition significantly. Notice that the notion of characteristic length scale becomes weaker as the parameter $\bar{\gamma}$ approaches the value $1$. 
\end{remark}
Concerning our notion of being mixed at scale $\lambda^k$, we simply require that the solution~$\rho$ has zero average on all cells $Q\in\mathcal{T}\!\!_{\lambda^k}$. We summarize our terminology in the following definition:
\begin{definition}[Mixed and un-mixed]
	\label{Def:MixedUnmixed}
	Let $\lambda>0$ be a tiling parameter and $\alpha\in(0,1)$ a universal constant. We say that a function~$\rho$ is both \textit{mixed} and \textit{un-mixed} at scale $\lambda^k$ if
	\begin{itemize}
		\item[(i)] $\int_{Q}\rho(x)\,dx=0$ for all $Q\in \mathcal{T}\!\!_{\lambda^k}$ and
		\item[(ii)] On every $Q\in\mathcal{T}\!\!_{\lambda^k}$, the function $\rho$ has the form
		\begin{equation*}
		\label{woreport}
		\rho|_Q=\begin{cases}1 & \mbox{on }A_Q\\ -1 & \mbox{on }A_Q^c\end{cases}
		\end{equation*} 
		where $LS\!_{Q}(A_Q)\geq \alpha\lambda^k$.
	\end{itemize}
\end{definition}

The following elementary lemma shows that, if a function is mixed at scale $\lambda$ according to Definition~\ref{Def:MixedUnmixed}(i), then both the geometric and the functional mixing scale are smaller than~$C\lambda^k$. See the Appendix of~\cite{CellP} for its proof.
 
 \begin{lemma}
	\label{tilingmixingPalen}
	Let $\rho$ be a bounded function such that 
	\begin{equation}
	\label{dunixPalen}
	\dashint_Q \rho\,dy=0
	\end{equation}
	for every tile $Q\in \mathcal{T}\!\!_{\lambda}$. Then there exist constants $C_1=C_1(\kappa)$ and $C_2=C_2(\|\rho\|_\infty)$ such that
	\begin{itemize}
		\item[(i)]
		$\mathcal{G}(\rho)\leq C_1\lambda$,
		\item[(ii)]
		$\|\rho\|_{\dot{H}^{-1}}\leq C_2 \lambda$.
	\end{itemize}
\end{lemma}

In the following definition we summarize our notion of generalized cellular structure.
\begin{definition}[Generalized cellular structure]
	\label{celltype}
	We say that a pair of a velocity field $u$ with corresponding solution $\rho$ to the continuity equation $\eqref{Cauchyprob}$ is of \textit{$(\lambda,\alpha,\bar{\gamma}, s, p, (T_n)_n, (\sigma(n))_n)$-generalized cellular type}, or simply of \textit{generalized cellular type} when parameters are fixed, if on any time interval $T_n\leq t < T_{n+1}$ and for any cell $Q\in \mathcal{T}\!\!_{\lambda^n}$ we can write $u$ and $\rho$ as
	\begin{equation}
	\label{TilVelSol2}
	u(t,x)=\frac{\lambda^n}{T_{n+1}-T_n}u_{Q,n}\left(\frac{t-T_n}{T_{n+1}-T_n},\frac{x-r_Q}{\lambda^n}\right)\hspace{0.5cm}\textnormal{and}\hspace{0.5cm}\rho(x,t)=\rho_{Q,n}\left(\frac{t-T_n}{T_{n+1}-T_n},\frac{x-r_Q}{\lambda^n}\right),
	\end{equation}
	where $r_Q$ is the center of the tile $Q$ and $(u_{Q,n},\rho_{Q,n})$ is a solution to \eqref{Cauchyprob} on $[0,1]\times \mathcal{Q}$, such that
	\begin{itemize}
		\item [(i)] $u_{Q,n}\in L_t^{\infty}(\dot{W}_x^{s,p})$ on the time interval $0\leq t\leq 1$ and is divergence free, $u_0=0$ on $\partial \mathcal{Q}$ and zero outside of $\mathcal{Q}$,
		\item [(ii)] $\rho_{Q,n}(0,\cdot)$ is both \textit{mixed} and \textit{un-mixed} at scale $\lambda^{\sigma(n)}$.
	\end{itemize}
\end{definition}
\begin{remark}\label{r:compa} Note that Definition~\ref{celltype} implies the following facts.
	\begin{itemize}
		\item [(i)] At time~$T_n$, the solution~$\rho$ is both \textit{mixed} and \textit{un-mixed} at scale~$\lambda^{n+\sigma(n)}$ (see Figure~\ref{yuumi}).
		\begin{figure}[H]
			\begin{center}
				\scalebox{0.3}{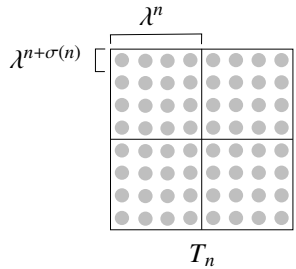}
			\end{center}
			\caption{Example of $\rho$ at iteration step $n$. \label{yuumi}}
		\end{figure}
		\item[(ii)] Since $(u,\rho)$ solves \eqref{Cauchyprob} and exploiting \eqref{TilVelSol2}, it follows that $\rho_{Q,n}(1,\cdot)$ is both \textit{mixed} and \textit{un-mixed} at scale $\lambda^{1+\sigma(n+1)}$.
		\item[(iii)] Definition~\ref{celltype} is a generalization of the cellular structure in~\cite{CellP}, which corresponds to the case $\sigma(n)=0$ for all $n\in\mathbb{N}$. 
	\end{itemize}
\end{remark}

\section{Statement and proof of the Main Theorem}
\label{s:main}

In this section we state and prove our main theorem, which asserts that under a fixed palenstrophy constraint it is not possible to mix exponentially fast by using a generalized cellular structure.

\begin{mainthm}
	Let $(u,\rho)$ be a solution of generalized cellular type (as in Definition~\ref{celltype}) of the continuity equation~\eqref{Cauchyprob} such that $\|u(t,\cdot)\|_{\dot{W}^{s,p}}$ is bounded uniformly in time, for some $s,p>1$. Then the decay of both the geometric and the functional mixing scale is slower than exponential.
\end{mainthm}

\begin{proof} We divide the proof in two steps.

\medskip

\noindent\textbf{Minimal Cost Estimate.} Let us fix $n\in\mathbb{N}$ and a cell $Q\in\mathcal{T}\!\!_{\lambda^n}$. Let $k\in\mathbb{N}$ to be chosen later. Note that by Definition~\ref{celltype},  by concatenating on the time intervals $[T_n,T_{n+1}]$, $[T_{n+1},T_{n+2}]$, \ldots, $[T_{n+k-1},T_{n+k}]$ solutions as in~\eqref{TilVelSol2}, on $[T_n,T_{n+k}]\times  Q$ the solution $(u,\rho)$ has the form
\begin{equation*}
\label{TilVelSol3}
u(t,x)=\frac{\lambda^n}{T_{n+k}-T_n}u_{Q,n,k}\left(\frac{t-T_n}{T_{n+k}-T_n},\frac{x-r_Q}{\lambda^n}\right)\hspace{0.5cm}\textnormal{and}\hspace{0.5cm}\rho(x,t)=\rho_{Q,n,k}\left(\frac{t-T_n}{T_{n+k}-T_n},\frac{x-r_Q}{\lambda^n}\right),
\end{equation*}
where $\rho_{Q,n,k}(0,\cdot):\mathcal{Q}\to\mathbb{R}$ is both mixed and un-mixed at scale~$\lambda^{\sigma(n)}$ and $\rho_{Q,n,k}(1,\cdot)$ is both mixed and un-mixed at scale~$\lambda^{k+\sigma(n+k)}$. We will prove the following minimal cost estimate
\begin{equation}
\label{MinCostEst}
c_1\log\left(\bar{c_2}\frac{\mix(\rho_{Q,n,k}(0,\cdot))}{\mix(\rho_{Q,n,k}(1,\cdot))}\right)\leq c_1\log\left(\frac{c_2}{\lambda^{k+\sigma(n+k)-\sigma(n)}}\right)=M_{n,k}\leq \int_{0}^{1}\|\nabla u_{Q,n,k}(t,\cdot)\|_{L^p(\mathcal{Q})}\,dt
\end{equation}
for $k$ large enough. The main ingredient for this step is the regularity result by Crippa and De Lellis~\cite{LellisCrippa} in Theorem~\ref{Thm:CrippaDeLellisRegOG}.

Let $A=\{x \text{ such that } \rho_{Q,n,k}(0,x)=1\}$. Since $\rho_{Q,n,k}(0,\cdot)$ is un-mixed at scale~$\lambda^{\sigma(n)}$, for every $Q\in\mathcal{T}\!\!_{\lambda^{\sigma(n)}}$ there exists a ball $B_Q=B_Q(x,r)\subset Q$, such that $r\geq\omega\lambda^{\sigma(n)}\alpha$ (where $\omega\in(0,1)$ will be chosen later  depending on $\bar\gamma$) and 
\begin{equation}
\label{fbfli}
\frac{|A\cap B_Q|}{|B_Q|}>1-\frac{1-\kappa}{2}\cdot\bar{\gamma}\geq 1-\frac{\bar{\gamma}}{2}\,,
\end{equation}	
with $\alpha$ the constant in Definition~\ref{Def:MixedUnmixed}. For every $Q\in\mathcal{T}\!\!_{\lambda^{\sigma(n)}}$ we consider the ball $\tilde{B}_Q=B(x,\sigma)\subset B_Q$ with the same center as $B_Q$ and such that $|\tilde{B}_Q|=\big((3+\bar{\gamma})/4\big)|B_Q|$. Thus 
\begin{equation}
\label{Nullsechs}
|A\cap\tilde{B}_Q|\geq 
|A \cap B_Q| - |B_Q \setminus \tilde{B}_Q| \geq
\left(\frac{3-\bar\gamma}{4}\right) |B_Q| =
%|\tilde{B}_Q|-\frac{\bar{\gamma}}{2}|B_Q|\geq \left(\frac{3+\bar{\gamma}}{4}-\frac{\bar{\gamma}}{2}\right)\pi r^2=
\left(\frac{3-\bar{\gamma}}{4}\right)\pi r^2\geq \left(\frac{3-\bar{\gamma}}{4}\right)\pi \omega^2\lambda^{2\sigma(n)}\alpha^2  .
\end{equation}
We now define 
\begin{equation}
B=\bigcup\limits_{Q\in\mathcal{T}\!\!_{\lambda^{\sigma(n)}}}B_Q\hspace{0.5cm}\text{ and }\hspace{0.5cm}\tilde{B}=\bigcup_{Q\in\mathcal{T}\!\!_{\lambda^{\sigma(n)}}}\tilde{B}_Q \,.
\end{equation}
Since the balls $B_Q$ are pairwise disjoint, by \eqref{fbfli} it follows that
\begin{equation}
\label{WhyG}
\frac{|A\cap B|}{|B|}>1-\frac{1-\kappa}{2}\cdot\bar{\gamma}\geq 1-\frac{\bar{\gamma}}{2}\,,
\end{equation}
and by \eqref{Nullsechs} we have that
\begin{equation}
\label{Sonate}
|A\cap\tilde{B}|\geq \left(\frac{3-\bar{\gamma}}{4}\right)\pi \omega^2\alpha^2\,.
\end{equation}
Furthermore, note that there exists a constant $1>C=C(\bar{\gamma})>0$ (in fact, $C(\bar\gamma)=1-\sqrt{( 3+\bar{\gamma})/4}$) such that
\begin{equation}
\label{panzamensch}
\dist(A\cap\tilde{B},B^c)\geq C(\bar{\gamma})r \geq \omega C(\bar{\gamma})\lambda^{\sigma(n)}\alpha\, ,
\end{equation}
which we need later on in the proof. 

We now focus on $\rho_{Q,n,k}$ at time $t=1$. We set $\Phi=X(1,\cdot)$, where $X:[0,1]\times \mathcal{Q}\to \mathcal{Q}$ is the  flow associated to $u_{Q,n,k}$. For notational purposes, we denote by $\delta=\lambda^{k+\sigma(n+k)}$. Let $\left\lbrace Q_i\right\rbrace_{i=1}^{1/\delta^2}$ denote the sub-squares in the tiling $\mathcal{T}\!\!_{\delta}$. For $i=1,\ldots, 1 / \delta^2$, denote
\begin{equation*}
A_i:=\Phi\left(A\cap\tilde{B}\right)\cap Q_i\, .
\end{equation*} 
Note that the $A_i$'s are disjoint for different $i$'s and by~\eqref{Sonate} we have that
\begin{equation}
\label{SeventySix}
\sum_{i=1}^{1/\delta^2}|A_i|=|A\cap\tilde{B}|\geq\left(\frac{3-\bar{\gamma}}{4}\right) \omega^2\alpha^2\pi\,.
\end{equation}
Let $G_i :=\Phi(A^c)\cap Q_i$. We further decompose each $G_i$ into $G_i^{\text{ good}}$ and $G_i^{\text{ bad}}$, where
\begin{equation*}
G_i^{\text{ good}}:=\Phi(A^c)\cap Q_i\cap\Phi(B^c) \hspace{0.5cm}\text{ and }\hspace{0.5cm}G_i^{\text{ bad}}:=\Phi(A^c)\cap Q_i\cap\Phi(B).
\end{equation*}
Recall that \eqref{WhyG} yields that
\begin{equation}
\label{LeftUnsaid}
\sum_{i=1}^{1/\delta^2}|G_i^{\text{bad}}|=|A^c\cap B|\leq\frac{\bar{\gamma}}{2}\alpha^2\pi\,,
\end{equation}
which we will use later in the proof. 

By Theorem~\ref{Thm:CrippaDeLellisRegOG}, for any $\eta>0$, there exists $E\subset \mathcal{Q}$ with $|E|<\eta$, such that
\begin{equation}
\label{ovaso}
\Lip(\Phi^{-1}|_{E^c})\leq \exp\left(\frac{C}{\eta^{1/p}}\int_{0}^{1}\|\nabla u_{Q,n,k}(s,\cdot)\|_{L^p}\,ds\right)\, .
\end{equation}
Let us fix
\begin{equation}
\label{Eta}
\eta=(1-\bar{\gamma})\frac{3}{16}\alpha^2\pi
\end{equation}
and let $E$ be the corresponding set as in~\eqref{ovaso}. For such $E$, we claim that there exists some index $i\in \{1,\ldots ,1/\delta^2\}$, such that \emph{both} $A_i\setminus E$ and $G_i^{\text{ good}}\setminus E$ are nonempty. Once this is proved, taking any points $\tilde{x}\in A_i\setminus E$ and $y\in G_i^{\text{ good}}\setminus E$ yields
\begin{equation*}
\Phi^{-1}(\tilde{x})\in A\cap \tilde{B}\hspace{0.5cm}\text{ and }\hspace{0.5cm}\Phi^{-1}(y)\in A^c\cap B^c \,.
\end{equation*}
Note that by \eqref{panzamensch} this implies that $|\Phi^{-1}(\tilde{x})-\Phi^{-1}(y)|\geq C\lambda^{\sigma(n)}\alpha$ (where $C$  depends on $\bar{\gamma}$ and $\omega$). On the other hand, since $\tilde{x},y\in Q_i$ we have that $|\tilde{x}-y|\leq\sqrt{2}\delta=\sqrt{2}\lambda^{k+\sigma(n+k)}$. Hence, by Theorem~\ref{Thm:CrippaDeLellisRegOG}, this yields
\begin{equation*}
\frac{C\alpha}{\lambda^{k+\sigma(n+k)-\sigma(n)}}\leq\frac{|\Phi^{-1}(\tilde{x})-\Phi^{-1}(y)|}{|\tilde{x}-y|}\leq\Lip(\Phi^{-1}|_{E^c})\leq \exp\left(\frac{C}{\eta^{1/p}}\int_{0}^{1}\|\nabla u_{Q,n,k}(s,\cdot)\|_{L^p}\,ds\right)\, .
\end{equation*}
Choosing $k$ large enough, such that the left hand side is larger than one, by the above inequality and by our particular choice of $\eta$ in~\eqref{Eta} there exist positive constants $c_1,c_2$ depending on $\bar{\gamma}$ and $\alpha$ only (which are in turn fixed parameters) such that
\begin{equation}
\label{ihro}
M_{n,k}=c_1\log\left(\frac{c_2}{\lambda^{k+\sigma(n+k)-\sigma(n)}}\right)\leq\int_{0}^{1}\|\nabla u_{Q,n,k}(s,\cdot)\|_{L^p}\,ds\,.
\end{equation}
Noting that,  by Remark~\ref{r:compa}(i), there exist constants $C_1$ and $C_2$ such that
\begin{equation}
\mix(\rho_{Q,n,k}(0,\cdot))\leq C_1\lambda^{\sigma(n)}\hspace{1cm}\text{ and }\hspace{1cm}\mix\rho_{Q,n,k}(1,\cdot))\geq C_2\lambda^{k+\sigma(n+k)}\,,
\end{equation} 
by \eqref{ihro} and setting $\bar{c_2}= c_2 C_2 / C_1$ we obtain the minimal cost estimate~\eqref{MinCostEst}. 

\smallskip

It only remains to prove the claim. For each $i=1,\ldots,1/\delta^2$, note that
\begin{equation}
\label{EightyOne}
\min\big\{|A_i\setminus E|,\, |G_i^{\text{ good}}\setminus E|\big\}
\geq \min\big\{ |A_i|, |G_i|\big\}-|G_i^{\text{ bad}}|-|Q_i\cap E| 
 =  |A_i | - |G_i^{\text{ bad}}| - |Q_i \cap E| \,.
\end{equation}
Summing \eqref{EightyOne} up for $i=1,\ldots,1/\delta^2$, by \eqref{SeventySix}, \eqref{LeftUnsaid} and \eqref{Eta} and setting $\omega^2=(1+\bar{\gamma})/2$ this yields
\begin{equation*}
\begin{split}
\sum_{i=1}^{1/ \delta^2}\min\big\{|A_i\setminus E|,\, |G_i^{\text{ good}}\setminus E|\big\}&\geq |A\cap\tilde{B}|-|A^c\cap B|-|E|
\geq\left[\left(\frac{3-\bar{\gamma}}{4}\right) \omega^2-\frac{\bar{\gamma}}{2}-(1-\bar{\gamma})\frac{3}{16}\right] \alpha^2\pi \\
&\geq\left[\frac{3}{4}\omega^2-\frac{3\bar{\gamma}}{4}-(1-\bar{\gamma})\frac{3}{16}\right] \alpha^2\pi\ 
=\left[1-\bar{\gamma}\right]\frac{3}{16}\alpha^2\pi>0 \,,
\end{split}
\end{equation*}
hence the claim is proved.

\medskip

\noindent\textbf{Scaling Computation.} We compute $\|\nabla_x^s u(t,\cdot)\|_{L^p(\mathcal{Q})}$ on the time interval $[T_n,T_{n+k})$, where $k=k(n)$ is as in the previous step. For notational convenience we denote $|[T_n,T_{n+k}]|=\tau_{n,k}$. Remember that for any tile $Q\in\mathcal{T}\!\!_{\lambda^n}$ the velocity field $u$ has the form
\begin{equation*}
u(t,x)|_Q\equiv\frac{\lambda^n}{\tau_{n,k}}u_{Q,n,k}\left(\frac{t-T_n}{\tau_{n,k}},\frac{x-r_Q}{\lambda^n}\right)\, .
\end{equation*}
 We will perform the following step for the case~$s\in\mathbb{N}$. The computation for the case~$s\in\mathbb{R}\setminus\mathbb{N}$ is analogous, however due to the non-locality of the fractional derivative the first equality in~\eqref{bigcalc} changes to "$\geq$". For more details, we refer to~\cite{CellP}, where this step is performed in detail for the case $s\in\mathbb{R}\setminus\mathbb{N}$. Using the change of variables, the Poincar\'e inequality, and the Jensen inequality we get
\begin{equation}
\label{bigcalc}
\begin{split}
\int\limits_{T_n}^{T_{n+k}}\|\nabla^s u(t,\cdot) \|_{L^p(\mathcal{Q})}^p\,dt&=\int\limits_{T_n}^{T_{n+k}}\sum_{Q\in\mathcal{T}\!\!_{\lambda^n}}\int\limits_{Q}\left| \nabla_x^s\left[\frac{\lambda^n}{\tau_{n,k}}u_{Q,n,k}\left(\frac{t-T_n}{\tau_{n,k}},\frac{x-r_Q}{\lambda^n}\right)\right]\right|^p\,dx\,dt\\
&=\frac{1}{\lambda^{(s-1)np}}\frac{\tau_{n,k}}{\tau_{n,k}^p}\sum_{Q\in\mathcal{T}\!\!_{\lambda^n}}\int\limits_{0}^{1}\int\limits_{Q}\left| \left(\nabla_x^s u_{Q,n,k} \right)\left(t,\frac{x-r_Q}{\lambda^n}\right)\right|^p\,dx\,dt\\
&=\frac{\tau_{n,k}}{(\lambda^{s-1})^{np}\tau_{n,k}^p}\lambda^{2n}\sum_{Q\in\mathcal{T}\!\!_{\lambda^n}}\int\limits_{0}^{1}\int\limits_{\mathcal{Q}}\left| \left(\nabla_x^s  u_{Q,n,k} \right)\left(t,x\right)\right|^p\,dx\,dt\,\\
&\stackrel{P}{\geq}\frac{C_{p,s}\tau_{n,k}}{(\lambda^{s-1})^{np}\tau_{n,k}^p}\lambda^{2n}\sum_{Q\in\mathcal{T}\!\!_{\lambda^n}}\int\limits_{0}^{1}\int\limits_{\mathcal{Q}}\left| \left(\nabla_x  u_{Q,n,k} \right)\left(t,x\right)\right|^p\,dx\,dt\,\\
&\stackrel{J}{\geq} \frac{C_{p,s}\tau_{n,k}}{(\lambda^{s-1})^{np}\tau_{n,k}^p}M_{n,k}^p\, ,
\end{split}
\end{equation}
where $M_{n,k}$ is as in \eqref{MinCostEst}, relative to the basic building blocks $u_{Q,n,k}$ used on $[T_n,T_{n+k}]$ . Now since $|[T_n,T_{n+1}]|=\tau_{n,k}$, together with \eqref{bigcalc} we get that
\begin{equation}
\label{endlichdeshiat}
\left\|\|\nabla_x^s u(t,\cdot)\|_{L^p(\mathcal{Q})}^p\right\|_{L_t^{\infty}(T_n,T_{n+1})}\geq \left(\frac{1}{\lambda^{s-1}}\right)^{np}\frac{1}{\tau_{n,k}^p} C_{p,s} M_{n,k}^p \,.
\end{equation}
Taking the $p$-th root of the above inequality, using the minimal cost estimate~\eqref{ihro} and the fact that the $W^{s,p}$ norm is uniformly bounded in time, this implies that
\begin{equation}
\label{Keedi}
C\geq \left(\frac{1}{\lambda^{s-1}}\right)^{n}\frac{1}{\tau_{n,k}} C_{p,s} c_1\log\left(\frac{c_2}{\lambda^{k+\sigma(n+k)-\sigma(n)}}\right)\,.
\end{equation}
Note that there exist constants $C_1$ and $C_2$ such that
\begin{equation}
\mix(\rho(T_n,\cdot))\leq C_1\lambda^{n+\sigma(n)}\hspace{1cm}\text{ and }\hspace{1cm}\mix(\rho(T_{n+k},\cdot))\geq C_2\lambda^{n+k+\sigma(n+k)}\,,
\end{equation}
which implies
\begin{equation}
\label{helleen}
\frac{\mix(\rho(T_n,\cdot))}{\mix(\rho(T_{n+k},\cdot))}\leq C\frac{\lambda^{n+\sigma(n)}}{\lambda^{n+k+\sigma(n+k)}}=\frac{C}{\lambda^{k+\sigma(n+k)-\sigma(n)}}\,.
\end{equation}
Using~\eqref{helleen} in~\eqref{Keedi} and denoting $\beta=\lambda^{s-1}\in (0,1)$ this yields
\begin{equation}
\label{zuul}
\exp[C\beta^n(T_{n+k}-T_n)]\geq C\frac{\mix(\rho(T_n,\cdot))}{\mix(\rho(T_{n+k},\cdot))}\,. 
\end{equation}
By Theorem~\ref{Thm:CrippaDeLellisRegOG} and the fact that the velocity field~$u$ is bounded in $W^{1,p}$ uniformly in time we have
\begin{equation}
\label{lowerzuul}
\mix(\rho(T_n))\geq \gamma_1\exp(-\alpha_1T_n)
\end{equation}
for all $n\in\mathbb{N}$. Now let us assume that the mixing scale of $\rho$ decays exponentially in time, i.e.~that there exists $\alpha_2>0$, such that
\begin{equation}
\label{upperzuul}
\mix(\rho(T_n))\leq \gamma_2\exp(-\alpha_2T_n)
\end{equation}
for all $n\in\mathbb{N}$. By~\eqref{lowerzuul} we must have  $\alpha_1\geq\alpha_2$. Then, by~\eqref{zuul}, \eqref{lowerzuul} and \eqref{upperzuul} it follows that
\begin{equation}
\begin{split}
\gamma_2\exp(-\alpha_2T_{n+k})&\geq\mix(\rho(T_{n+k},\cdot)) \\
&\geq C\mix(\rho(T_n,\cdot))\exp[-C\beta^n(T_{n+k}-T_n)]\\
&\geq C\exp[-C\beta^n(T_{n+k}-T_n)] \, \gamma_1\exp(-\alpha_1T_n)\\
&=C\gamma_1\exp \big[-\big(C\beta^n(T_{n+k}-T_n)+\alpha_1T_n \big)\big] 
\end{split}
\end{equation}
and hence
\begin{equation}
\label{ching}
\bar{C}\geq\exp\left[(\alpha_2-C\beta^n)T_{n+k}-(\alpha_1-C \beta^n)T_n\right]
\end{equation}
for a positive constant $\bar{C}$. This is however a contradiction: if we choose $n$ large enough in such a way that both $(\alpha_2-C\beta^n)>0$ and $(\alpha_1-C \beta^n)>0$, we can then let $k$ go to infinity (which in turn lets $T_{n+k}$ go to infinity) and we then see that the right hand side of~\eqref{ching} goes to infinity.
\end{proof}

\end{document}